\documentclass{amsart}

\usepackage{amssymb,graphicx}

\newtheorem{thm}{Theorem}[section]

\newtheorem{lem}[thm]{Lemma}
\newtheorem{cor}[thm]{Corollary}

\begin{document}

\title[A Fast Algorithm to the Conjugacy Problem]
{A Fast Algorithm to the Conjugacy Problem\\ on Generic Braids}

\author[K.~H.~KO]{Ki Hyoung KO}
\author[J.~W.~LEE]{Jang Won LEE}
\email{\{knot,leejw\}@knot.kaist.ac.kr}
\address{Department of Mathematics, Korea Advanced Institute
of Science and Technology, Daejeon, 305-701, Korea}
\maketitle

\begin{abstract}
Random braids that are formed by multiplying randomly chosen
permutation braids are studied by analyzing their behavior under
Garside's weighted decomposition and cycling. Using this analysis,
we propose a polynomial-time algorithm to the conjugacy problem that
is successful for random braids in overwhelming probability. As
either the braid index or the number of permutation-braid factors
increases, the success probability converges to 1 and so, contrary
to the common belief, the distribution of hard instances for the
conjugacy problem is getting sparser. We also prove a conjecture by
Birman and Gonz\'{a}lez-Meneses that any pseudo-Anosov braid can be
made to have a special weighted decomposition after taking power and
cycling. Moreover we give polynomial upper bounds for the power and
the number of iterated cyclings required.
\end{abstract}

\section{Preliminaries and introduction}

Recently the braid groups have become a potential source for
cryptography, especially, for public-key cryptosystems (see
\cite{AAG,KLC} for few). The braid groups have two important features
that are useful for cryptography. Each word can be quickly put into
a unique canonical form, which provides a fast algorithm not only
for the word problem but also for the group operation (see
\cite{G,T,EM}). On the other hand no polynomial-time solution to the
conjugacy problem in the braid group is known, which provides many
interesting one-way functions for public-key cryptosystems.

Before we discuss the history and the main result, we quickly
introduce the terminologies and basic facts about braid groups.
Artin who first studied braids systematically in the early 20th
century proved that the group $B_n$ of $n$-strand braids can be
given by the following presentation:
$$ B_n =\left\langle \sigma_1, \cdots, \sigma_{n-1} \biggm|
\begin{array}{ll}
 \sigma_j\sigma_i =\sigma_i\sigma_j &{\rm if}~~ \vert i-j \vert > 1 \\
  \sigma_i\sigma_j\sigma_i = \sigma_j\sigma_i\sigma_j &{\rm if}~~ \vert i-j
  \vert = 1\
  \end{array}
  \right\rangle.$$
The monoid given by the same presentation is denoted by $B_n^+$
whose elements will be called {\em positive braids}.

A partial order $\prec$ on $B_n^+$ can be given by saying $x\prec y$
for $x, y\in B_n^+$ if $x$ is a ({\em left}) {\em subword} of $y$,
that is, $xz=y$ for some $z\in B_n^+$. Given $x,y\in B_n^+$, the
({\em left}) {\em join} $x\vee y$ of $x$ and $y$ is the minimal
element with respect to $\prec$ among all $z$'s satisfying that
$x\prec z$ and $y\prec z$, and the ({\em left}) {\em meet} $x\wedge
y$ of $x$ and $y$ is the maximal element with respect to $\prec$
among all $z$'s satisfying that $z\prec x$ and $z\prec y$. Even
though ``left" is our default choice, we sometimes need the
corresponding right versions: the partial order $\prec_R$ of being a
{\em right subword}, the {\em right join} $\vee_R$, and the {\em
right meet} $\wedge_R$. For example, $x\prec_R y$ if $zx=y$ for some
$z\in B_n^+$.

The {\em fundamental braid} $\Delta=(\sigma_1\cdots
\sigma_{n-1})(\sigma_1\cdots \sigma_{n-2})\cdots
(\sigma_1\sigma_2)\sigma_1$ plays an important role in the study of
$B_n$. Since it represents a half twist as a geometric braid,
$x\Delta =\Delta\tau(x)$ for any braid $x$ where $\tau$ denotes the
involution of $B_n$ sending $\sigma_i$ to $\sigma_{n-i}$. It also
has the property that $\sigma_i\prec\Delta$ for each
$i=1,\ldots,n-1$. Since the symmetric group $\Sigma_n$ is obtained
from $B_n$ by adding the relations $\sigma_i^2=1$, there is a
quotient homomorphism $q:B_n\to \Sigma_n$. For $S_n=\{x\in B_n^+\mid
x\prec\Delta\}$, the restriction $q:S_n\to\Sigma_n$ becomes a 1:1
correspondence and an element in $S_n$ is called a {\em permutation
braid}.

A product $ab$ of a permutation braid $a$ and a positive braid $b$
is ({\em left}) {\em weighted}, written $a\lceil b$, if $a^* \wedge
b=e$ where $e$ denotes the empty word and $a^*=a^{-1}\Delta$ is the
{\em right complement} of $a$. Each braid $x\in B_n$ can be uniquely
written as
$$x=\Delta^{u}x_1x_2\cdots x_{k}$$
where for each $i=1,\ldots,k$, $x_i\in S_n\setminus
\{e,\Delta\}$ and $x_i\lceil x_{i+1}$. This decomposition is called
the ({\em left}) {\em weighted form} of $x$ \cite{G,T,EM}. Sometimes
the first and the last factors in a weighted form are called the
{\em head} and the {\em tail}, denoted by $H(x)$ and $T(x)$,
respectively. The weighted form provides a solution to the word
problem in $B_n$ and the integers $u$, $u+k$ and $k$ are
well-defined and are called the {\em infimum}, the {\em supremum} and
the {\em canonical length} of $x$, denoted by $\inf(x)$, $\sup(x)$
and $\ell(x)$, respectively.

Given $x=\Delta^{u}x_1x_2\cdots x_{k}$ in its weighted form, there
are two useful conjugations of $x$ called the {\em cycling} $\mathbf
c(x)$ and the {\em decycling} $\mathbf d(x)$ defined as follows:
$${\mathbf c}(x)=\Delta^{u}x_2\cdots x_{k}\tau^{u}(x_1)=
\tau^{u}(H(x)^{-1})x\tau^{u}(H(x)),$$
$${\mathbf d}(x)=\Delta^{u}\tau^{u}(x_{k})x_1\cdots x_{k-1}=
T(x)xT(x)^{-1}.$$

A braid $x$ is {\em cyclically weighted} if its weighted form
$x=\Delta^{u}x_1x_2\cdots x_{k}$ has the property that
$x_{k}\lceil\tau^{u}(x_1)$. A braid $x$ is {\em weakly
cyclically weighted} if $H(x\tau^u(x_1))=x_1=H(x)$. A cyclically
weighted braid is clearly weakly cyclically weighted. When
$\ell(x)=1$, the two properties are equivalent and they require
$x_1\lceil\tau^{u}(x_1)$.

Let $\inf_c(x)$ and $\sup_c(x)$ respectively denote the maximum of
infimums and the minimum of supremums of all braids in the conjugacy
class $C(x)$ of $x$. A typical solution to the conjugacy problem in
the braid group $B_n$ is to generate a finite set uniquely
determined by a conjugacy class. Historically, the following four
finite subsets of the conjugacy class $C(x)$ of $x\in B_n$ have been
used in this purpose:

The {\em summit set}
$$SS(x)=\{y\in C(x)\mid \inf(y)=\mbox{$\inf_c$}(x)\}$$
was used by Garside in \cite{G} to solve the conjugacy problem in
$B_n$ for the first time. The {\em super summit set}
$$SSS(x)=\{y\in C(x)\mid \inf(y)=\mbox{$\inf_c$}(x)\mbox{ and }\sup(y)=\mbox{$\sup_c$}(x)\}$$
was used by El-Rifai and Morton in \cite{EM} to improve Garside's
solution. The {\em reduced super summit set}
$$RSSS(x)=\{y\in C(x)\mid \mathbf c^M(y)=y=\mathbf d^N(y)
\mbox{ for some positive integers }M,N\}$$ was used by Lee in his
Ph.D. thesis \cite{L} to give a polynomial-time solution to the
conjugacy problem in $B_4$. Finally the {\em ultra summit set}
$$USS(x)=\{y\in SSS(x)\mid \mathbf c^M(y)=y
\mbox{ for some positive integer }M\}$$ was used by Gebhardt in
\cite{Ge} to propose a new algorithm together with experimental data
demonstrating the efficiency of his algorithm. Clearly
$$RSSS(x)\subset USS(x)\subset SSS(x)\subset SS(x),$$
and $RSSS(x)=USS(x)$ if $x$ is cyclically weighted.
 All four
invariant sets enjoy the property that if $a^{-1}ya\in P$ and
$b^{-1}yb\in P$ for $y\in P$ and $a,b\in S_n$ then $(a\wedge
b)^{-1}y(a\wedge b)\in P$ where $P$ denotes one of invariant sets.
So for $y\in P$, there is a minimal element $a\in S_n$ such that
$a^{-1}ya\in P$. Franco and Gonz\'{a}lez-Meneses\cite{FG} first
proved this property for the super summit set  and then
Gebhardt\cite{Ge} did it for the ultra summit set. Using this
property, they were able to generate an invariant set more
efficiently. Unfortunately there is no estimate for the sizes of the
invariant sets and so we do not know the complexity of any algorithm
based on the generation of an invariant set.

In this paper we survey a fast algorithm to the conjugacy problem on
generic braids. In the algorithmic sense, generic braids means
random braids that are built by multiplying randomly chosen
permutation braids. In the dynamical sense, generic braids means
pseudo-Anosov braids. In Section 2, we first give a combinatorial analysis
on random braids to find out how quickly the head of a random braid
becomes stable as the braid index or the canonical length increases.
Then we show that a random braid is cyclically weighted up to cycling in
an overwhelming probability so that its $RSSS(x)$ is predictable and
small. Using this, we propose a polynomial-time algorithm to the
conjugacy problem for random braids. Some of proofs are omitted or brief
in this section and full proofs will appear elsewhere.

In Section 3, we show that some
power of a pseudo-Anosov braid is always cyclically weighted up to
cycling and we also give upper bounds for the necessary exponent
and the necessary number of iterated cyclings. Our upper bounds are
polynomial in canonical length so that there would be a
polynomial-time solution to the conjugacy problem for pseudo-Anosov
braids once the size of reduced super summit sets are known to be
polynomial in canonical length. Finally we give an example of a
cyclically weighted pseudo-Anosov braid whose reduced super summit
set is relatively large to show there are still some more work required
to give a good estimate of the size of reduced super summit sets.

\section{A fast algorithm to the conjugacy problem for random braids}

We assume that the permutation braids in $S_n$ are uniformly distributed
so that each permutation braid can be chosen with an equal probability
of $1/n!$. We consider random braids that are formed by multiplying $k$
factors, each of which is a permutation braid chosen randomly from
$S_n$. Random braids need not be positive and nonpositive random braids
are obtained by multiplying a random (negative) power of the fundamental
braid $\Delta$. Since $\Delta$
commutes with any braid up to the involution $\tau$, a power of
$\Delta$ can be ignored in most of the discussions so that random braids
are assumed to be positive. In this section, we study the behavior of
random braids with respect to two parameters $k$ and $n$. We reveal
some unexpected facts regarding the braid index $n$.

For integers $1\le i< j \le n$, we say that a positive $n$-braid $x$
begins with an inversion $(i,j)$ if the head $H(x)$ exchanges $i$ and
$j$ as a permutation. For permutation braids $x_1, x_2,\ldots,x_k$
chosen randomly from $S_n$, let $D(n,k,(i,j))$ denote the
probability that $x_1x_2\cdots x_k$ begins with the inversion
$(i,j)$. In particular an inversion $(i,i+1)$ that a positive braid
$x$ begins with is called a {\em descent} of $x$ and $\mathcal D(x)$
denotes the set of all descents of $x$. We will write
$D(n,k,(i,i+1)):=D(n,k,i)$. Then $D(n,k):=\sum_{i=1}^{n-1}D(n,k,i)$
denotes the average number of descents of a random braid
$x_1x_2\cdots x_k$. It is easy to see that $D(n,1,(i,j))=1/2$ for
$i=1,\ldots, n-1$ since the $i$-th and $j$-th strings in $x_1$ cross
each other in the probability $1/2$. Thus $D(n,1)=(n-1)/2$. We need
more delicate combinatorial analysis to obtain an estimate of
$D(n,k)$ for $k\ge 2$ that is sharp enough to be useful. In fact we
will give an estimate on how fast $d(n,k):=D(n,k)-D(n,k-1)$, the
average contribution to descents of the product $x_1x_2\cdots x_k$
by the last factor $x_k$, approaches to $0$ as either $n$ or $k$
increases.

\begin{lem}\label{d(n,2)}We have
$$
D(n,2,i)=\frac12+\frac{1}{n(n-1)}\sum_{k=0}^{n-2}\frac{n-k-1}{\binom{n-2}{k}}
\sum_{j=0}^{k}\frac{\binom{i-1}{j}\binom{n-i-1}{k-j}}{\binom{k+2}{j+1}}$$
and $\displaystyle d(n,2)\le
\frac{1}{n}\sum_{k=0}^{n-2}\frac{n-k-1}{k+2}$. In particular,
$d(n,2)$ is in $O(\log n)$ and is not a decreasing function.
\end{lem}
\begin{proof}
For a braid $x\in B_n$, $\hat x:\{1,2,\ldots,n\}\to\{1,2,\ldots,n\}$
denotes the permutation $q(x)\in \Sigma_n$. In order that $\sigma_i$
is a descent of $x=x_1x_2$ contributed by $x_2$, all the
following two conditions must hold.
\begin{enumerate}
\item[(i)]$\hat x_1(i)<\hat x_1(i+1)$;
\item[(ii)]If $\hat x_1^{-1}(a)\le i$ and $\hat x_1^{-1}(b)\ge i+1$
for $a$, $b$ with $\hat x_1(i)\le a,b\le \hat x_1(i+1)$, then $\hat
x_2(b)<\hat x_2(a)$.
\end{enumerate}
The condition (i) contributes 1/2. The number of choices for $\hat
x_1(i)$ and $\hat x_1(i+1)$ can be expressed in two distinct ways:
$\binom{n}{2}=\sum_{k=0}^{n-2}(n-k-1)$ where $k=\hat x_1(i+1)-\hat
x_1(i)-1$. The number of choices for $k$ integers sent between $\hat
x_1(i)$ and $\hat x_1(i+1)$ by $\hat x_1$ can be also expressed in
two ways:
$\binom{n-2}{k}=\sum_{j=0}^{k}\binom{i-1}{j}\binom{n-i-1}{k-j}$.
Then the $k+2$ integers $\hat x_1(i),\hat x_1(i)+1,\ldots, \hat
x_1(i+1)$ are divided into two groups such that the first group
consists of $j+1$ integers whose preimage under $\hat x_1$ is less
than or equal to $i$, and the remaining $k-j+1$ integers have
preimages greater than or equal to $i+1$. The condition (ii) requires
that $\hat x_2$ permutes the $k+2$ numbers so that each image of
the first group is larger than all images of the
second group. The claimed formula for $D(n,2,i)$ in Lemma should
be clear now. The rest of proof is technical and omitted.
\end{proof}

In general, we have the following properties that are extremely
useful to give an estimate for an upper bound of $D(n,k,i)$.
\begin{lem}\label{inequalities}
\begin{enumerate}
\item $D(n,k,i)=D(n,k,n-i)$ and $D(n,k,i)\ge D(n,k,j)$ for $1\le i< j \le \lfloor n/2\rfloor$
\item $D(n,k,(i,j))\le D(n+i-j+1,k,i)$ (The equality holds for $k=1,2$)
\end{enumerate}
\end{lem}

\begin{proof} Since  $\mathcal D(x_1\ldots x_k)=\mathcal D(x_1H(x_2\ldots x_k))$,
the argument for random braids of two factors in the previous lemma
can similarly be applied to show (1). For a random $n$-braid
$x=x_1\ldots x_k$ made of $k$ permutation braids, let $x'$ be the
$(n+i-j+1)$-braid obtained from $x$ by deleting $j-i-1$ strings from
the $(i+1)$-th to $(j-1)$-th. If $x$ begins with a inversion $(i,j)$,
then $x'$ must have the descent $\sigma_i$. The converse is also
true for $k=1, 2$. This proves (2).
\end{proof}

\begin{thm}\label{recursive bound}
For all $n\ge 2$, $k\ge 2$, and $1\le i\le n-1$, $D(n,k,i)$ is
recursively bounded above by $$\frac12 + \frac
2{n(n-1)}\sum_{a=0}^{n-2}\frac{n-a-1}{(a+2)!}D(n-a,k-1,1)\sum_{b=0}^{a}(a-b+1)!b!D(n-a-1,k-2,1)^b$$
where $D(n,0,1)=0$ and $D(n,1,1)=1/2$ for all $n$.
\end{thm}
\begin{proof}
Again since $\mathcal D(x_1\ldots x_k)=\mathcal D(x_1H(x_2\ldots
x_k))$, a typical usage of induction on $k$ together with
inequalities in Lemma~\ref{inequalities} gives a proof. The details
are omitted.
\end{proof}

As a corollary, we have the following estimate for $d(n,3)$. This is
rather surprising because the total number of descents of a random
$n$-braid contributed by the third factor (and by all following
factors) eventually decreases to $0$ as the braid index increases.
The maximum occurs at $n=9$ and this means that 9-braids are the
most well-mixed in their weighted forms, for example, when two
braids are multiplied.

\begin{cor} We have
$$D(n,3,1)\le\frac12 + \frac{\ln n}{n-1} + \frac{3(\ln n)^2}{n(n-1)}$$
and so asymptotically $\displaystyle d(n,3)\le \frac{3(\ln
n)^2}{n}$.
\end{cor}
\begin{proof} Omitted
\end{proof}

Even though a recursive upper bound is given in
Theorem~\ref{recursive bound}, it is difficult to describe an upper
bound for $D(n,k,i)$ as a neat formula.  Instead we use
$(n-1)(D(n,k,1)-D(n,k-1,1))$ as an estimate for an upper bound of
$d(k,n)$ and present a table for these upper bounds for some choices of
$(n,k)$ that are relevant to Gebhardt's experiment in \cite{Ge}. The
table shows that $d(n,k)$ converge to $0$ as $k$ increases and
moreover the larger the $n$ becomes the faster it converges to $0$.

Since $d(n,k)$ quickly converges to $0$ and $D(n,k)$ is much less
than $n-1$, it is extremely difficult to produce $\Delta$
by multiplying randomly chosen permutation $n$-braids unless the number
of chosen permutation $n$-braids is comparable to $n!$. Thus we assume
in the rest of the article that $\inf(x_1\cdots x_k)=0$.
Since $\sup(x)=-\inf(x^{-1})$, we may assume that $\sup(x_1\cdots x_k)=k$
as well.

\begin{table*}[h]\label{table}
\begin{center}
\begin{tabular}{|c|c|c|c|c|c|}
%
  \hline
  $k \backslash n$ & 4 & 6 & 8 & 10&15 \\
  \hline
  2 & $6.04\times 10^{-1}$ & $8.58\times 10^{-1}$ & $1.06$ & $1.22$ & $1.54$ \\
  5 & $9.08\times 10^{-2}$ & $1.67\times 10^{-1}$ & $2.01\times 10^{-1}$ & $2.01\times 10^{-1}$ &$1.57\times 10^{-1}$ \\
  10 & $3.00 \times 10^{-3}$ & $7.17\times 10^{-3}$ & $1.19\times 10^{-2}$ & $1.57\times 10^{-1}$ & $1.54\times 10^{-2}$ \\
  20 & $2.91\times 10^{-6}$ & $7.03\times 10^{-6}$ & $1.21\times 10^{-5}$ & $1.70\times 10^{-5}$ & $2.18\times 10^{-5}$ \\
  30 & $2.85\times 10^{-9}$ & $6.86\times 10^{-9}$ & $1.18\times 10^{-8}$ & $1.66\times 10^{-8}$ & $2.12\times 10^{-8}$ \\
  40 & $2.78\times 10^{-12}$ & $6.70\times 10^{-12}$ & $1.16\times 10^{-11}$ & $1.62\times 10^{-11}$ & $2.07\times 10^{-11}$ \\
  50 & $3.00\times 10^{-15}$ & $6.11\times 10^{-15}$ & $1.24\times 10^{-14}$ & $1.50\times 10^{-14}$ & $2.02\times 10^{-14}$ \\
  \hline
  \hline
  $k \backslash n$ & 20&30&50&75&100 \\
  \hline
  2 &  $1.78$ & $2.13$ &$2.59$ &$2.97$&$3.24$\\
  5 &  $1.18\times 10^{-1}$ & $7.4\times 10^{-2}$ &$3.82\times 10^{-2}$ & $2.17\times 10^{-2}$& $1.43\times 10^{-2}$\\
  10 &  $9.33\times 10^{-3}$ & $3.21\times 10^{-3}$ & $8.61\times 10^{-4}$& $3.22\times 10^{-4}$&$1.65\times 10^{-4}$\\
  20 &  $1.70\times 10^{-5}$ & $6.85\times 10^{-6}$ &$1.74\times 10^{-6}$ & $6.25\times 10^{-7}$&$3.14\times 10^{-7}$\\
  30 &  $1.66\times 10^{-8}$ & $6.77\times 10^{-9}$ &$1.73\times 10^{-9}$& $6.20\times 10^{-10}$ &$3.11\times 10^{-10}$\\
  40 &  $1.62\times 10^{-11}$ & $6.61\times 10^{-12}$ &$1.67\times 10^{-12}$ & $6.08\times 10^{-13}$&$2.97\times 10^{-13}$\\
  50 &  $1.48\times 10^{-14}$ & $6.44\times 10^{-15}$ & $<10^{-15}$& $<10^{-15}$& $<10^{-15}$\\
  \hline
\end{tabular}
\end{center}
\caption{Upper bounds for $d(n,k)$}
\end{table*}

We now observe some of the properties that random braids enjoy with an
overwhelming probabilities.
We will use the notation $\mathrm{Prob}[S(x):x]$ or simply
$\mathrm{Prob}[S(x)]$ to denote the probability that the statement
$S(x)$ is true for a random choice of $x$.

\begin{lem} \label{stability}For randomly chosen
permutation $n$-braids $x_1,\ldots, x_k$, let $x=x_1\cdots x_k$.
Then the probability that the following equivalent properties hold
is greater than $1-\min\{d(n,k),1\}$:
\begin{enumerate}
\item For a
randomly chosen permutation $n$-braid $a$,
$H(x)=H(xa)$;\\
\item For a randomly chosen permutation $n$-braid $a$,
$T(x)=T(ax)$.
\end{enumerate}
\end{lem}

\begin{proof}
If $H(x)\not=H(xa)$, then $\mathcal D(x_2\ldots x_k)\subsetneqq
\mathcal D(x_2\ldots x_ka)$.
$$\mathrm{Prob}[\mathcal D(x_2\ldots x_k)\not=
\mathcal D(x_2\ldots x_ka):(x_2,\ldots,x_k,a)]\le \min\{d(n,k),
1\}$$ because $d(n,k)$ is the average contribution to descents by
the $k$-th factor which is $a$. Thus (1) follows.

Consider $(ax)^{-1}=x_k^*\tau(x_{k-1}^*)\cdots
\tau^{k-1}(x_1^*)\tau^{k}(a^*)\Delta^{-(k+1)}$. Then
$$T(ax)^*=H(x_k^*\tau(x_{k-1}^*)\cdots
\tau^{k-1}(x_1^*)\tau^{k}(a^*))$$ and similarly
$T(x)^*=H(x_k^*\tau(x_{k-1}^*)\cdots \tau^{k-1}(x_1^*))$. If
$x_1,\ldots, x_k, a$ are random, so are $x_k^*,
\tau(x_{k-1}^*),\ldots, \tau^{k-1}(x_1^*),\tau^{k}(a^*)$. Thus (2)
follows since $T(ax)=T(x)$ if and only if
$H(x_k^*\tau(x_{k-1}^*)\cdots
\tau^{k-1}(x_1^*)\tau^{k}(a^*))=H(x_k^*\tau(x_{k-1}^*)\cdots
\tau^{k-1}(x_1^*))$.
\end{proof}

\begin{lem}\label{random and WCW}For $k\ge 3$ and randomly chosen permutation $n$-braids
$x_1,\ldots, x_k$  and a randomly chosen integer $u$, the
probability that $x=\Delta^{u}x_1\cdots x_k$ is weakly cyclically
weighted is greater than $1-2d(n,k)$. In particular
$\mathrm{Prob}[x\in SS(x):x]\ge 1-2d(n,k)$
\end{lem}

\begin{proof} For the simplicity of notation, we assume $u=0$.
Lemma~\ref{stability} implies
\begin{align*}
\mathrm{Prob}&[H(x_1\cdots x_k)=H(x_1\cdots x_k x_1\cdots
x_k):(x_1,\ldots,x_k)]\\
&\ge 1-\min\{d(n,k)+\cdots+d(n,2k),1\}
\end{align*}


According to our estimate via Theorem~\ref{recursive bound},
$1>d(n,k)\ge 2d(n,k+1)$ for $k\ge 3$. Thus
$$1-\min\{d(n,k)+\ldots+d(n,2k),1\}\ge 1-2d(n,k).$$
If $x$ is weakly cyclically weighted,
$\inf(x)=\inf(\mathbf{c}^i(x))$ for all $i>0$ and so $x\in SS(x)$.
\end{proof}

\begin{cor}\label{random and SSS}For $k\ge 3$ and randomly chosen permutation $n$-braids
$x_1,\ldots, x_k$ and a randomly chosen integer $u$, the probability
that $x=\Delta^{u}x_1\cdots x_k\in SSS(x)$ is greater than
$1-4d(n,k)$.
\end{cor}

\begin{proof}
By Lemma \ref{random and WCW},
\begin{align*}
\mathrm{Prob}[x\notin SSS(x)]&\leq\mathrm{Prob}[x\notin
SS(x)\mbox{~or~} x^{-1}\notin
SS(x^{-1})]\\
&\leq\mathrm{Prob}[x\notin SS(x))+\mathrm{Prob}(x^{-1}\notin
SS(x^{-1})]\\
&=2d(n,k)+2d(n,k)
\end{align*}
and so
$$\mathrm{Prob}[x\in SSS(x)]>1-4d(n,k).$$
\end{proof}

\begin{lem} \label{cut-head}For $k\ge 3$ and randomly chosen
permutation $n$-braids $x_1,\ldots, x_k$, let $x=x_1\cdots x_k$.
Then the probability that the following equivalent properties hold
is greater than $1-2d(n,k)$:
\begin{enumerate}
\item For any permutation braid $a$, $H(x)=H(xa)$ or $\tau(H(xa))$;\\
\item For any permutation braid $a$, $T(x)=T(ax)$.
\end{enumerate}
\end{lem}
\begin{proof} The argument is similar to Lemma~\ref{stability} but
the difference is that we need to assume that the probability that a
pair of strands has a crossing in $a$ is 1 for (1) and 0 for (2). On
the other hand, it was 1/2 in Lemma~\ref{stability}. Then this lemma
becomes obvious.
\end{proof}

\begin{thm}\label{random and CW}For $k\ge 12$ and randomly chosen permutation $n$-braids
$x_1,\ldots, x_k$and a random integer $u$, the probability that
${\bf{c}}^j(\Delta^{u}x_1\cdots x_k)$ is cyclically weighted for
some $j\leq \lfloor \frac{k}{2}\rfloor$ is greater than
$1-2d(n,\lfloor \frac{k}{4}\rfloor)$.
\end{thm}

\begin{proof} For the simplicity, we again assume $u=0$. We also assume
that $x\in SSS(x)$ and this happens with the probabiltiy $\ge
1-4d(n,k)$. Let $y_1\cdots y_k$ be the weighted form of $x=x_1\cdots
x_k$, $h=\lfloor\frac{k}{2}\rfloor$, and $q=\lfloor\frac{k}{4}\rfloor$.

Set $a=T(x_1\cdots x_{h})^*\wedge H(x_{h+1}\cdots x_{k})$. Then we have
$$H(y_{h+1}\cdots y_ky_1\cdots y_{h})=H(y_{h+1}\cdots y_ky_1)=
H(a^{-1}x_{h+1}\cdots x_{k}y_1)$$ and
$$y_{h+1}=H(y_{h+1}\cdots y_{k})=H(a^{-1}x_{h+1}\cdots x_{k}).$$
By Lemma~\ref{cut-head},
$$\mathrm{Prob}[H(x_{h+q+1}\cdots x_k)=
H(x_{h+q+1}\cdots x_ky_1)]\ge 1-2d(n,k-h-q)$$
and
$$\mathrm{Prob}[T(a^{-1}x_{h+1}\cdots x_{h+q})=
T(a^{-1}x_{h+q}\cdots x_{h+q}]\ge 1-2d(n,q).$$
So
$$\mathrm{Prob}[H(a^{-1}x_{h+1}\cdots x_k)=H(a^{-1}x_{h+1}\cdots x_ky_1)]
\ge 1-2d(n,q).$$
Thus
$$\mathrm{Prob}[H(y_{h+1}\cdots y_ky_1\cdots y_{h})=y_{h+1}]
\ge 1-2d(n,q).$$


Similarly, we have
$$\mathrm{Prob}[T(y_{h+1}\cdots y_ky_1\cdots y_{h})=y_{h}]
\ge 1-2d(n,q).$$

Since $y_{j}\lceil y_{j+1}$, the probability that
${\bf{c}}^j(\Delta^{u}x_1\cdots x_k)$ is cyclically weighted is
greater than or equal to $1-2d(n,\lfloor \frac{k}{4}\rfloor)$. We note
that $1>2d(n,\lfloor \frac{k}{4}\rfloor)>4d(n,k)$ for $k\ge 12$ and so
our assumption $x\in SSS(x)$ makes no dfference.
\end{proof}

We now know from Theorem~\ref{random and CW} that a random braid can
be made cyclically weighted by a small number of iterated cyclings
with an overwhelming probability. A cyclically weighted braid $x$
already belongs to $USS(x)=RSSS(x)$ and hence the conjugacy problem can be
solved by generating $USS(x)$. In the remaining of this section, we will show $USS(x)$
is very small for a random braid $x$, in fact $|USS(x)|\le
2\ell(x)$, with an overwhelming probability.

Let $P$ denote one of the conjugacy invariant sets $SS, SSS, USS,
RSSS$ and let $y\in P(x)$. If a nontrivial positive $n$-braid
$\gamma$ satisfies $\gamma^{-1}y\gamma \in P(x)$, $\gamma$ is
called a {\em $P$-conjugator} of $y$. A $P$-conjugator $\gamma$ of
$y$ is {\em minimal} if either $\gamma\prec \beta$ or
$\gamma\wedge \beta=e$ for each positive braid $\beta$ with
$\beta^{-1}y\beta \in P(x)$. In fact it is not hard to see that a
minimal $P$-conjugator satisfies $\gamma\prec
\tau^{\inf(y)}(H(y))$ or $\gamma\prec T(y)^*$ or both (For
example, see \cite{FG}). A conjugator $\gamma$ satisfying
$\gamma\prec \tau^{\inf(y)}(H(y))$ (or $\gamma\prec T(y)^*$,
respectively) will be called a {\em cut-head} (or {\em add-tail})
conjugator. In particular, if $y$ is cyclically weighted and
$\gamma$ is its $P$-conjugator then it can not be both cut-head
and add-tail since $T(y)\lceil \tau^{\inf(y)}(H(y))$, that is,
$T(y)^*\wedge \tau^{\inf(y)}(H(y))=e$. If $\gamma$ is a
$USS$-conjugator of a cyclically weighted braid $y$, it is also a
$RSSS$-conjugator and $\gamma^{-1}y\gamma$ is cyclically weighted.
We note that if $\gamma$ is an add-tail conjugator of $y$, then
$\gamma$ is a cut-head conjugator of $y^{-1}$.

\begin{thm} \label{random and orbit number1} For $k\ge 3$ and randomly chosen permutation $n$-braids
$x_1,\ldots, x_k$ and a randomly chosen integer $u$, assume that
$y\in USS(\Delta^{u}x_1\cdots x_k)$ is cyclically weighted and $t$
is a USS-minimal cut-head (or add-tail, respectively) conjugator of
$y$. Then the probability that $t=\tau^u(H(y))$ (or $t=T(y)$) is
greater than $1-2d(n,k-1)$.
\end{thm}

\begin{proof} If $T(t^{-1}y)=T(y)$ and $t\precneqq\tau^u(H(y))$,
$t^{-1}yt$ can not be in $SSS(\Delta^{u}x_1\cdots x_k)$ since
$\sup(t^{-1}yt)=\sup(y)+1$. Now the conclusion is immediate from
Lemma~\ref{cut-head}.
\end{proof}

\begin{cor} \label{random and orbit number3} For $k\ge 12$ and randomly
chosen permutation $n$-braids $x_1,\ldots, x_k$ and an integer $u$,
let $x=\Delta^{u}x_1\cdots x_k$. Then the probability that
$\mathrm{USS}(x)=\mathrm{O}(y)\cup\mathrm{O}(\tau(y))$ for some
$y\in USS(x)$ is greater than
$1-2d(n,\lfloor\frac{k}{4}\rfloor)$ where
$\mathrm{O}(y)=\{\mathbf{c}^i(y)\mid i>0\}$ is a finite set called
the {\em cycling orbit} of $y\in USS(x)$.
\end{cor}
\begin{proof} Immediate from Theorem~\ref{random and CW} and
Theorem~\ref{random and orbit number1} since $1>2d(n,\lfloor\frac{k}{4}\rfloor)>
2d(n,k-1)$ for $k\ge 12$.
\end{proof}

Given a
random braid $x\in B_n$, an algorithm to generate $USS(x)$ is now extremely simple.
In fact, it proceeds as follows:
\begin{enumerate}
\item Compute $y=\mathbf{c}^{j}(x)$ where
$j=\lfloor\frac{\ell(x)}{2}\rfloor$.
\item Output either
$\{\mathbf{c}^i(y),~\tau(\mathbf{c}^i(y))\mid 0\le i\le
\ell(x)-1\}$ if $\inf(y)$ is even, or
$\{\mathbf{c}^i(y)\mid 0\le i\le 2\ell(x)-1\}$ if $\inf(y)$ is odd.
\end{enumerate}

Since the operation to build a left canonical form has running time
$\mathcal{O}(k^2 n \log n)$ (see \cite{T}) and this algorithm
requires $k$ cycling operations, the overall running time is
$\mathcal{O}(k^3 n \log n)$. For a given $n$-braid $x$ of $k$ random
permutations, it generates $USS(x)$ successfully with probability
greater than $1-2d(n,\lfloor\frac{k}{4}\rfloor)$ by
Corollary~\ref{random and orbit number3}.

\section{Conjugacy problem for pseudo-Anosov braids}

As far as Garside's weighted decomposition is concerned, we will
show that pseudo-Anosov braids behave similarly to random braids
that we discussed. This is rather surprising because the dynamical
notion of generic braids are seemingly far from the combinatorial
notion. On the other hand, this may be natural in the sense that a
braid chosen randomly as a mapping class should be expected to be
pseudo-Anosov. J.~Gonz\'{a}lez-Meneses discovered a surprising
phenomenon that some power of any pseudo-Anosov braid is cyclically
weighted up to cyclings and J.~Birman announced this phenomenon as a
conjecture at the first East Asian School of Knot Theory and Related
Topics in 2004. We verify this conjecture and give upper
bounds for the exponent and the number of iterated cyclings
required. Recently the proposers independently
verified their conjecture in \cite{BGG}.
In short, we will prove that for any pseudo-Anosov braid
$x$, there are integers $1\le M\le D^3$ and $1\le N\le D^4\ell(x)$
such that $\mathbf c^N(x^M)$ is cyclically weighted where
$D=\frac{n(n-1)}{2}$. In \cite{KL}, We give a polynomial-time algorithm
to decide the dynamical type of any given braid by using this special
property and these bounds.

\begin{lem}\label{lem:wcwpower}Let $x$ be an $n$-braid. Then
there exists $y\in C(x)$ and an integer $1\le M\le D^3$ such that
$y^M$ is weakly cyclically weighted.
\end{lem}
\begin{proof}
It was shown in \cite{LL} that for any $n$ braid $x$, there exists
$y\in C(x)$ and an integer $M_1>0$ such that
$\inf((y^{M_1})^i)=i\inf(y^{M_1})$ for all $i\ge 1$. Let
$z=y^{M_1}$. Since any unexpected $\Delta$ can not be produced by
taking powers of $z$, there exists $0<M_2<D$ such that
$H(z^{M_2})=H(z^{k})$ for all $k\ge M_2$ as argued in \cite{BKL}.
Thus $z^{M_2}$ is weakly cyclically weighted.
\end{proof}

\begin{lem}\label{lem:CW}Let $x$ be an $n$-braid. If $x$ is cyclically
weighted and $\ell(x)\ge 2$, then every braid in $RSSS(x)$ is
cyclically weighted.
\end{lem}
\begin{proof}
Since $x$ is cyclically weighted, $x\in RSSS(x)$. It is enough to
show that $t^{-1}xt$ is cyclically weighted when $t$ is minimal
among conjugators such that $t^{-1}xt\in RSSS(x)$. Let
$x=\Delta^ux_1\cdots x_{\ell}$ be the weighted form that is
cyclically weighted. Then either $\tau^u(t)\prec x_1$ or $t\prec
x_{\ell}^*$ since $t$ is minimal. Suppose that $\tau^u(t)\prec x_1$.
Then
$$t^{-1}xt=\Delta^{u}(t_1^{-1}x_1t_2)\cdots
(t_{\ell}^{-1}x_{\ell}\tau^{-u}(t_1))$$ is the weighted form where
$t_1=\tau^u(t)$ and $t_i=x_{i}\wedge x_{i-1}^*\tau(t_{i-1})$ for $1<
i \le \ell$. If $t^{-1}xt$ is not cyclically weighted, then
$(t_{\ell}^{-1}x_{\ell}\tau^{-u}(t_1)) (\tau^{-u}(t_1^{-1}x_1t_2))$
is not weighted. Thus $$\mathbf
c(t^{-1}xt)=\Delta^u(t_2^{-1}x_2t_3)\cdots
(t_{\ell-1}^{-1}x_{\ell-1}t_{\ell})(t_{\ell}^{-1}x_{\ell}
\tau^{-u}(t_1'))\tau^{-u}(t_1'^{-1}x_1t_2)$$ is the weighted form
and $t_1\prec t_1'\prec x_1$ because
$$\tau^{-u}(t_1)^{-1}x_{\ell}^*\tau(t_{\ell})\wedge
tau^{-u}(t_1)^{-1}\tau^{-u}(x_1t_2)\ne e$$ and so
$\tau^{-u}(t_1)\precneqq x_{\ell}^*\tau(t_{\ell})\wedge
\tau^{-u}(x_1)=\tau^{-u}(t_1')$. Since the last $\ell-2$ factors of
$t^{-1}xtT(t^{-1}xt)^{-1}$ are equal to the first $\ell-2$ factors of
$\mathbf c(t^{-1}xt)$, $\mathbf
c^{\ell}(t^{-1}xt^{-1})=t_1'^{-1}xt_1'$ and $\tau^u(t)\precneqq
t_1'\prec x_1$. This implies that $\mathbf c^{2\ell
i}(t^{-1}xt^{-1})=s_{(i)}^{-1}xs_{(i)}$ for $i\ge 0$, where
$s_{(0)}=t$ and $s_{(i)}\precneqq s_{(i+1)} \prec \tau^{-u}(x_1)$.
Thus $s_{(j)}=\tau^{-u}(x_1)$ for some $j>0$ and so
$s_{(i)}=\tau^{-u}(x_1)$ for $i>j$. Since $t^{-1}xt\in RSSS(x)$ and
so $t^{-1}xt=\mathbf c^{2\ell j'}(t^{-1}xt)$ for some $j'>j$,
$t^{-1}xt=\mathbf c^{2\ell
j'}(t^{-1}xt)=\tau^{-u}(x_1)^{-1}x\tau^{-u}(x_1)=\mathbf c(x)$. But
this is a contradiction since $\mathbf c(x)$ is cyclically weighted. Thus
$t^{-1}xt$ must be cyclically weighted. If $t\prec x_{\ell}^*$ then
$\tau^{-u-\ell}(t)\prec H(x^{-1})$. Since $x^{-1}\in RSSS(x^{-1})$
and $x^{-1}$ is cyclically weighted, so is $t^{-1}x^{-1}t$ by the
above. Thus $t^{-1}xt$ is cyclically weighted.
\end{proof}

\begin{thm}[Birman-Gonz\'{a}lez-Meneses Conjecture]\label{lem:BG-conj} Let $x$ be
a pseudo-Anosov $n$-braid. Then there is a positive integer $M$ such
that every braid in $RSSS(x^M)$ is cyclically weighted.
\end{thm}
\begin{proof}
By Lemma \ref{lem:wcwpower}, there exists $y\in C(x)$ and an integer
$K>0$ such that $y^K$ is weakly cyclically weighted and
$\ell(y^K)\ge 2$. Since $y^K\in SS(x^K)$ and $SS(x^K)$ is finite,
$\mathbf c^{j}(y^K)=\mathbf c^{j+N}(y^K)$ for some $j, N>0$. Set
$z=\mathbf c^{j}(y^K)$ and $\mathbf c^i(z)=\Delta^uZ_iQ_i$ where
$Z_i$ and $Q_i$ are the head and the rest of $\mathbf c^i(z)$,
respectively. Since $z=\mathbf c^{N}(z)=\gamma^{-1}z\gamma$,
$\gamma$ is an element of the centralizer of $z$, where
$\gamma=\tau^{-u}(Z_0\cdots Z_{N-1})$. It is well-known that the
centralizer of $z$ is generated by a pseudo-Anosov $n$-braid
$\alpha$ and a periodic $n$-braid $\rho$ such that $\rho^p=\Delta^2$
(see \cite{GW}). Thus we can write $z=\rho^{u_1}\alpha^{v_1}$ and
$\gamma=\rho^{u_2}\alpha^{v_2}$. Then
$z^{v_2}=\rho^{u_1v_2-u_2v_1}\gamma^{v_1}$ and so
$z^{v_2p}=\Delta^{2(u_1v_2-u_2v_1)}\gamma^{v_1p}$. On the other
hand, $Z_i \lceil Z_{i+1}$ for $i>0$ since $y^K$ is weakly
cyclically weighted and so is $z$. Since $Z_{N-1} \lceil Z_N$ and
$Z_{N}=Z_0$, $Z_{N-1} \lceil Z_0$. Thus $\gamma$ is cyclically
weighted and so $z^{v_2p}$ is cyclically weighted. Since
$\ell(z^{v_2p})\ge 2$, every braid in $RSSS(x^{v_2pK})$ is
cyclically weighted by Lemma \ref{lem:CW}.
\end{proof}

\begin{lem}\label{lem:cwpower}
Suppose that $x$ is a braid such that $x\in SSS(x)$ and
$\inf(x^i)=i\inf(x)$, $\sup(x^i)=i\sup(x)$ for $i\ge 1$. If $x^N$ is
cyclically weighted for some $N\ge 1$ then $x$ itself is cyclically
weighted.
\end{lem}
\begin{proof} Under the hypotheses, neither new $\Delta$'s can be
formed nor factors can be merged by taking powers. Thus
$\tau^{u(N-1)}(H(x))\prec H(x^N)$ and $T(x)\succ_R T(x^N)$.
Consequently $T(x^N)\lceil\tau^{uN}(H(x^N))$ implies
$T(x)\lceil\tau^{u}(H(x))$.
\end{proof}

\begin{cor}\label{thm:pacw}
Let $x$ be a pseudo-Anosov braid in $B_n$. Then $x^M$ is conjugate
to a cyclically weighted braid for some $1\le M\le D^2$. Moreover,
every braid in $RSSS(x^{M'})$ is cyclically weighted for some $1\le
M'\le 2D^2$.
\end{cor}

\begin{proof}
In Theorem \ref{lem:BG-conj}, we have already proved the existence of such an $M$
and so we discuss the upper bound for $M$. By \cite{LL}, there
exists a positive integer $M\leq D^2$ and $y\in C(x)$ such that
$\inf((y^M)^i)=i\inf(y^M)$ and $\sup((y^M)^i)=i\sup(y^M)$ for $i\ge
1$. Let $z=y^M$. Since $z$ is pseudo-Anosov, $z^{M'}$ is conjugate
to a cyclically weighted braid for some $M'\ge 1$. Hence, $z$ is
conjugate to a cyclically weighted braid by Lemma \ref{lem:cwpower}.

On the other hand, every braid in $RSSS(z^2)$ is cyclically weighted
by Lemma \ref{lem:CW} since $\ell(z^2)\ge 2$. Thus $M'\le 2M$ and so
$1\le M'\le 2D^2$.
\end{proof}

The next theorem tells us how fast we can obtain a cyclically weighted 
braid that is conjugate to a power of a given pseudo-Anosov braid. In the
theorem, we assume that a braid is conjugate to a cyclically weighted braid 
instead of being pseudo-Anosov. This assumption is weaker because of  
Corollary~\ref{thm:pacw}. 

\begin{thm}\label{cycling bound and WCW}
Let $y$ be an $n$-braid such that $\inf((y)^i)=i\inf(y)$ and $\sup((y)^i)=i\sup(y)$
for all $i\ge 1$ and $y\in SSS(y)$ and let $x=y^{2D}$
If $x$ is conjugate to a cyclically weighted braid,
then a cyclically weighted braid must be obtained from $x$ by at most
$n!\,\ell(x)$ iterated cyclings.
\end{thm}
\begin{proof} It was proved in Lemma \ref{lem:CW} that
if $RSSS(x)$ contains at least one cyclically weighted braid, then
every braid in $RSSS(x)$ is cyclically weighted. Thus iterated
cyclings on $x$ must produce a cyclically weighted braid. Let
$y=\mathbf c^N(x)$ be the cyclically weighted braid obtained from
$x$ by the minimal number of iterated cyclings. Since $\inf(x)$ is even,
we assume $\inf(x)$ for  the sake of simplicity.

Let $y=y_1y_2\cdots y_k$ be the weighted form. First one can prove by
induction on $i\ge 1$ that for all $1\le i\le N$
$$\mathbf c^{N-i}(x)=a_{i}y_{[1-i]}z_{i}$$
for some permutation braid $a_{i}$ satisfying $y_{[2-k-i]}\cdots
y_{[-1-i]}y_{[-i]}\succ_R a_{i}\succneqq_R e$ and a positive braid
$z_{i}=y_{[2-k-i]}\cdots y_{[-1-i]}y_{[-i]}a_i^{-1}$ where $[m]$
denotes the integer between 1 and $k$ that equals $m$ mod $k$.
Then $\mathbf c^{N-i}(x)$ is completely determined by
choosing a nontrivial permutation braid $a_i$ satisfying
$$a_i\prec_R y_{[2-k-i]}\cdots y_{[-1-i]}y_{[-i]}.$$ For each
$1\le i\le \ell(x)$, there are at most $n!$ such choices. Thus
$N\le n!\,\ell(x)$. A complete proof appears in \cite{KL}. 
\end{proof}

Given any pseudo-Anosov braid, we now know that we are able to
generate a cyclically weighted braid that is conjugate to some power
of the given braid in polynomial time. Thus a polynomial-time
algorithm to solve the conjugacy problem for pseudo-Anosov braids
will be completed as soon as we know how to generate the whole set
$RSSS(x)$ for a pseudo-Anosov and cyclically-weighted braid $x$. If
$x$ is a pseudo-Anosov and cyclically-weighed braid obtained by
iterated cyclings on a product of randomly chosen permutation
braids, then $RSSS(x)$ has at most two cycling orbits in an
overwhelming probability. But there are plenty of pseudo-Anosov and
cyclically weighted braids whose reduced super summit set are not so
simple.

Consider the following permutation 7-braids:
\begin{equation*}
x_1=\sigma_2\sigma_1\sigma_3\sigma_2\sigma_5,~
x_2=\sigma_2\sigma_5\sigma_6,~ x_3=\sigma_2\sigma_6\sigma_5,~
x_4=\sigma_2\sigma_5\sigma_4\sigma_6\sigma_5.
\end{equation*}
\begin{figure}[h]
\center
\includegraphics{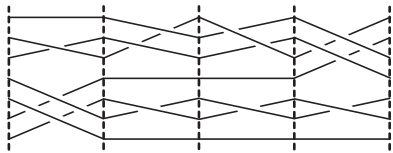}
\caption{\label{quasi-reducible} A pseudo-Anosov braid that is
``quasi-reducible"}
\end{figure}
Then $x=x_1x_2x_3x_4$ is pseudo-Anosov and cyclically weighted. But
$\sigma_5^{-1}x\sigma_5$ is again cyclically weighted and forms a
new cycling orbit in $RSSS(x)$. In fact, the number of cycling
orbits in $RSSS(x)$ is 10. We say that pseudo-Anosov braids of this kind
are {\em quasi-reducible} because they are almost reducible and
contain most of complications due to reducibility. Consequently we
still need more study to estimate the size of $RSSS(x)$ for a
pseudo-Anosov and cyclically weighted braid $x$.

\end{document}